\documentclass[12pt]{article}
\usepackage{amsfonts,amssymb}
\newcommand{\N}{{\mathbb N}}
\newcommand{\Z}{{\mathbb Z}}
\newcommand{\Fp}{{\mathbb F}_{\!p}}
\newcommand{\Q}{{\mathbb Q}}
\newcommand{\R}{{\mathbb R}}
\newcommand{\C}{{\mathbb C}}
\newcommand{\E}{{\cal E}}
\newcommand{\A}{{\cal A}}
\newcommand{\T}{\tilde{T}}
\newcommand{\ftilde}{f^*}
\newcommand{\ftildebarT}{\overline{f}^*_{\!\!T}}

\newcommand{\I}{\tilde{I}}
\newcommand{\abar}{\overline{a}}

\newcommand{\Jbar}{\overline{J}}
\newcommand{\Pbar}{\overline{P}}
\newcommand{\Upsilonbar}{\overline{\Upsilon}}
\newcommand{\fbar}{\overline{f}}
\newcommand{\gbar}{\overline{g}}
\renewcommand{\hbar}{\overline{h}}

\newcommand{\proof}{\noindent{\em Proof: }}
\newcommand{\ra}{\rightarrow}
\newcommand{\m}{\medskip}
\newcommand{\qed}{\hspace{\fill}$\square$}

\renewcommand{\theequation}{\thesection.\arabic{equation}}
\newtheorem{theorem}[equation]{Theorem}
\newtheorem{lemma}[equation]{Lemma}
\newtheorem{prop}[equation]{Proposition}
\newtheorem{cor}[equation]{Corollary}

\newenvironment{example}{\noindent\refstepcounter{equation}{\bf Example \theequation} }{}

\title{Signed Shape Tilings of Squares
}
\author{Kevin Keating\thanks{Partially
supported by NSF grant 9500982.} \\
Department of Mathematics \\
University of Florida \\
Gainesville, FL 32611 \\
USA \\
{\tt keating@math.ufl.edu}}
\date{}
\begin{document}

\maketitle

\begin{abstract}
Let $T$ be a tile made up of finitely
many rectangles
whose corners have rational coordinates
and whose sides are parallel to the
coordinate axes.  This paper gives
necessary and sufficient conditions for
a square to be tilable by finitely many
$\Q$-weighted tiles with the same shape
as $T$, and necessary and sufficient conditions for a
square to be tilable by finitely many
$\Z$-weighted tiles with the same shape
as $T$.  The main tool we use is a
variant of F.\,W.\,Barnes's algebraic
theory of brick packing, which
converts tiling problems into problems
in commutative algebra.
\end{abstract}

\section{Introduction} \label{intro}
\setcounter{equation}{0}

     In \cite{dehn} Dehn proved that an
$a\times b$ rectangle $R$ can be tiled
by finitely many nonoverlapping squares
if and only if $a/b$ is rational.
More generally, suppose we allow the squares to have
weights from $\Z$.  An arrangement of
weighted squares is a tiling of $R$
if the sum of the weights of
the squares covering a region is 1
inside of $R$ and 0 outside.
Dehn's argument applies in this more
general setting, and shows that $R$ has a
$\Z$-weighted tiling by squares if and
only if $a/b$ is rational.
In \cite{st} this result is generalized
to give necessary and sufficient
conditions for a rectangle $R$ to be tilable
by $\Z$-weighted rectangles with particular shapes.
In this paper we consider a related
question: Given a
tile $T$ in the plane made up of finitely
many weighted rectangles, is there a weighted
tiling of a square by tiles with the
same shape as $T$?

\begin{figure}
\begin{center}
\begin{picture}(400,150)
\put(35,0){\vector(0,1){150}}
\put(15,20){\vector(1,0){150}}
\put(255,0){\vector(0,1){150}}
\put(235,20){\vector(1,0){150}}
\put(45,35){\framebox(29.7,105){3}}
\put(75,105){\framebox(60,35){3}}
\put(100,65){\framebox(70,39.5){$-2$}}
\put(265,35){\framebox(30,69.5){3}}
\put(265,105){\framebox(90,35){3}}
\put(320,65){\framebox(70,39.5){$-2$}}
\end{picture}
\end{center}
\caption{Two rectangle decompositions of the same
$\Z$-weighted tile.}
\label{tile}
\end{figure}
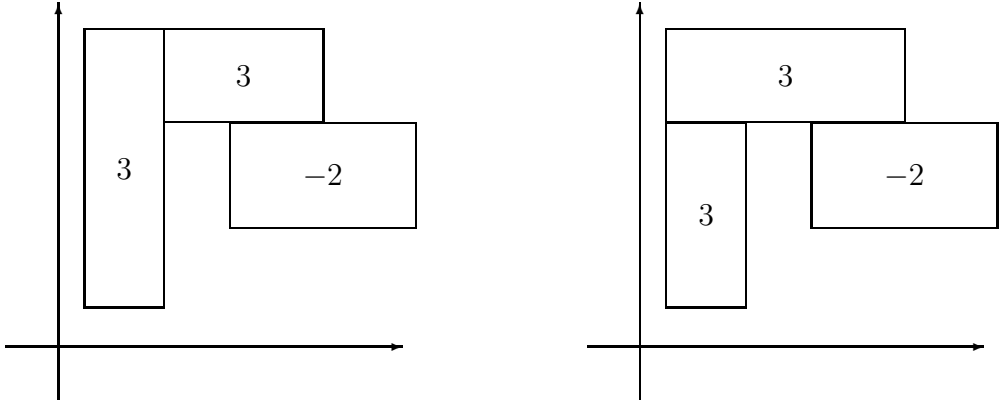

     We define a {\em rectangle} in $\R\times\R$
to be a product $[b_1,b_2)\times[c_1,c_2)$ of
half-open intervals, with $b_1<b_2$ and $c_1<c_2$.
Let $A$ be a commutative ring with unity.  An {\em
$A$-weighted tile} is represented by a finite
$A$-linear combination $L=a_1R_1+\cdots+a_nR_n$
of disjoint rectangles.  Associated to each such
$L$ there is a function $f_L:\R^2\ra A$ which is
supported on $\cup\,R_i$ and whose value on
$R_i$ is $a_i$.  We say that $L_1$ and
$L_2$ represent the same tile if
$f_{L_1}=f_{L_2}$.  An example of a
$\Z$-weighted tile is given in
Figure~\ref{tile}.
We may form the sum $T_1+T_2$ of two weighted tiles
$T_1,T_2$ by superposing them in the natural way.
For $a\in A$ the tile $aT$ is formed
from $T$ by multiplying all the weights
of $T$ by $a$.
The set of all $A$-weighted tiles forms
an $A$-module under these operations.

     Let $U$ be an $A$-weighted tile and
let $\{T_{\lambda}:\lambda\in\Lambda\}$ be
a set of $A$-weighted tiles.  We say
that the set
$\{T_{\lambda}:\lambda\in\Lambda\}$
{\em $A$-tiles} $U$ if there are
weights $a_1,\dots,a_n\in A$ and tiles
$\T_1,\dots,\T_n$, each of which is a
translation of some $T_{\lambda_i}$,
such that $a_1\T_1+\cdots+a_n\T_n=U$.
Note that we are allowed to use as
many translated copies of each
prototile $T_{\lambda}$ as we need, but
we are not allowed to rotate or reflect
the prototiles.
Given an $A$-weighted tile $T$ and
a real number $\rho>0$ 
we define $T(\rho)$ to be the image
of $T$ under the rescaling
$(x,y)\mapsto(\rho x,\rho y)$.  We
say that an $A$-weighted tile $T'$ has
the same shape as $T$ if there exists
$\rho>0$ such that $T'$ is a translation
of $T(\rho)$.  We say that $T$
$A$-{\em shapetiles} $U$ if
$\{T(\rho):\rho>0\}$ $A$-tiles $U$.  If
$U'$ has the same shape as $U$ then $T$
$A$-shapetiles $U'$ if and only if $T$
$A$-shapetiles $U$.

     In this paper we consider tiles $T$
constructed from rectangles whose
corners have rational coordinates.
We prove two main results about such
tiles.  First, we show that if $T$
is a $\Q$-weighted tile whose weighted
area is not 0, then $T$ $\Q$-shapetiles a
square.
Second, if $T$ is a $\Z$-weighted tile
we give necessary and sufficient conditions for $T$
to $\Z$-shapetile a square.

     The author would like to thank
Jonathan King for posing several
questions which led to this work.

\section{Polynomials and tiling}
\label{polynomials}
\setcounter{equation}{0}

     Say that $T$ is a {\em lattice tile}
if $T$ is an $A$-weighted tile made up
of unit squares in $\R^2$ whose corners
are in $\Z^2$.  We will associate a (generalized)
polynomial $f_T$
to each $A$-weighted lattice tile $T$.
Our approach is similar to that used by
F.\,W.\,Barnes in \cite{fwb1}, except that 
the polynomials that we construct differ
from Barnes's polynomials by a factor
$(X-1)(Y-1)$.  Including this extra
factor will allow us to generalize the
construction to non-lattice tiles at the
end of the section.

     Our polynomials will be elements of the ring
\begin{eqnarray*}
A[X^{\Z},Y^{\Z}]&:=&A[X,Y,X^{-1},Y^{-1}],
\end{eqnarray*}
which is naturally isomorphic to the
group ring of $\Z\times\Z$ with coefficients in
$A$.  To begin we associate the polynomial
$X^iY^j(X-1)(Y-1)$ to the unit square
$S_{ij}$ with lower
left corner $(i,j) \in \Z \times \Z$.
Given an $A$-weighted lattice tile
\begin{eqnarray*}
T&=&\sum_{i,j}w_{ij}S_{ij},
\end{eqnarray*}
by linearity we associate to $T$ the polynomial
\begin{eqnarray*}
f_T(X,Y)&=&\sum_{i,j}\,w_{ij}X^iY^j(X-1)(Y-1).
\end{eqnarray*}
One consequence of this definition is
that translating a tile by a vector ${(i,j)\in
\Z\times\Z}$ corresponds to multiplying its
polynomial by $X^iY^j$.  The map $T
\mapsto f_T$ gives an isomorphism
between the $A$-module of $A$-weighted
lattice tiles in the plane and the
principal ideal in $A[X^{\Z},Y^{\Z}]$
generated by $(X-1)(Y-1)$. \m

\begin{example} \label{rectangle}
Let $a,b,c,d$ be integers such that $a,b
\geq 1$ and let $T$ be
the $a \times b$ rectangle whose lower
left corner is at $(c,d)$.  Then the
polynomial associated to $T$ is
\begin{eqnarray*}
f_T(X,Y) &=&
\sum_{i=c}^{c+a-1}\;\sum_{j=d}^{d+b-1}\;X^{i}Y^{j}(X-1)(Y-1) \\
&=& X^cY^d(X^{a}-1)(Y^{b}-1).
\end{eqnarray*}
\vspace{-.7cm}
\end{example}

     In section~\ref{integer} we will
need to work with non-lattice tiles.  To
represent these more general tiles systematically
we introduce a new set
of building blocks to play the role that
the unit squares $S_{ij}$ play in the
theory of lattice tiles.
For $\alpha,\beta \in \R^{\times}$
let $R_{\alpha\beta}$ denote the oriented rectangle
with vertices $(0,0)$, $(\alpha,0)$,
$(\alpha,\beta)$, $(0,\beta)$.  Note
that if exactly $k$ of
$\alpha,\beta$ are negative then
$R_{\alpha\beta}$ is equal to $(-1)^k$ times
a translation of $R_{|\alpha|,|\beta|}$.
We can express any rectangle in terms of
the rectangles $R_{\alpha\beta}$: \m

\begin{example} \label{rectangle2}
Let $\alpha,\beta > 0$ and let $R_{\alpha\beta}'$
be the translation of the rectangle $R_{\alpha\beta}$
by the vector $(\sigma,\tau) \in \R^2$.  Then
$R_{\alpha\beta}'=R_{\alpha+\sigma,\beta+\tau} -
R_{\alpha+\sigma,\tau} - R_{\sigma,\beta+\tau} +
R_{\sigma\tau}$.  In particular, we have
$S_{ij}=R_{i+1,j+1}-R_{i+1,j}-R_{i,j+1}+R_{ij}$.
\end{example} \m

     In fact the following holds:

\begin{lemma} \label{welldef}
Every $A$-weighted tile $T$
can be expressed uniquely as an
$A$-linear combination of rectangles
$R_{\alpha\beta}$ with $\alpha,\beta \in
\R^{\times}$.
\end{lemma}

\proof By Example~\ref{rectangle2} every
rectangle is an $A$-linear combination
of the rectangles $R_{\alpha\beta}$.
Therefore every $A$-weighted tile is an
$A$-linear combination of the $R_{\alpha\beta}$.
Suppose
\begin{eqnarray*}
c_1R_{\alpha_1\beta_1}+c_2R_{\alpha_2\beta_2}+
\cdots+c_nR_{\alpha_n\beta_n} &=& 0
\end{eqnarray*}
is a linear relation such that
the pairs $(\alpha_i,\beta_i)$ are
distinct and
$c_i \not= 0$ for $1 \leq i \leq n$.
Choose $j$ to
maximize the distance from the origin to
the far corner $(\alpha_j,\beta_j)$ of
$R_{\alpha_j\beta_j}$.  None of the
other rectangles in the sum can overlap the region
around $(\alpha_j,\beta_j)$.
Since $c_j \not= 0$,
this gives a contradiction.
Therefore the set $\{R_{\alpha\beta} :
\alpha,\beta \in \R^{\times}\}$
is linearly independent over $A$, which implies the
uniqueness part of the lemma.~\qed \m

     In order to represent arbitrary
$A$-weighted tiles algebraically
we introduce a generalization of the
polynomials $f_T$.  Let $A[X^{\R},Y^{\R}]$ denote
the set of ``polynomials'' with coefficients from
$A$ where the exponents of $X$ and $Y$ are
allowed to be arbitrary real numbers.  The
natural operations of
addition and multiplication make
$A[X^{\R},Y^{\R}]$ a commutative
ring with unity.  The ring
$A[X^{\R},Y^{\R}]$ is naturally isomorphic to the
group ring of $\R\times\R$ with coefficients in
$A$, and contains $A[X^{\Z},Y^{\Z}]$ as a subring.

     For $\alpha,\beta \in \R^{\times}$
define $f_{R_{\alpha\beta}} =
(X^{\alpha}-1)(Y^{\beta}-1) \in A[X^{\R},Y^{\R}]$.
By Lemma~\ref{welldef} this definition
extends linearly to give a well-defined
element $f_T \in A[X^{\R},Y^{\R}]$ associated to
any $A$-weighted tile $T$.  It follows from
Example~\ref{rectangle2} that this definition agrees
with that given earlier if $T=S_{ij}$ is a unit lattice
square, and hence also if $T$ is any lattice tile.
The map $T \mapsto
f_T$ gives an isomorphism between the
$A$-module of $A$-weighted tiles
and an $A$-submodule of $A[X^{\R},Y^{\R}]$.  The next
lemma implies that this $A$-submodule is actually
an ideal in $A[X^{\R},Y^{\R}]$.

\begin{lemma} \label{translate}
Let $T$ be an $A$-weighted tile
and let $T'$ be the
translation of $T$ by the vector
$(\sigma,\tau) \in \R\times\R$.
Then $f_{T'} = X^{\sigma}Y^{\tau}f_T$.
\end{lemma}

\proof Let $R_{\alpha\beta}'$
be the translation of $R_{\alpha\beta}$ by
$(\sigma,\tau)$.  Using
Example~\ref{rectangle2} we get
\[f_{R_{\alpha\beta}'} \;\;=\;\;
X^{\sigma}Y^{\tau}(X^{\alpha}-1)(Y^{\beta}-1)
\;\;=\;\; X^{\sigma}Y^{\tau}f_{R_{\alpha\beta}},\]
so the lemma
holds for $T=R_{\alpha\beta}$.  Therefore by
Lemma~\ref{welldef} the lemma holds for all tiles
$T$.~\qed \m

     The next result gives a further relation
between ideals and tiling.

\begin{prop} \label{ideal}
Let $U$ be a tile, let
$\{T_{\lambda}:\lambda\in\Lambda\}$
be a collection of tiles,
and let $\I\subset A[X^{\R},Y^{\R}]$
be the ideal generated by the set
$\{f_{T_{\lambda}}:\lambda\in\Lambda\}$.
Then $\{T_{\lambda}:\lambda\in\Lambda\}$
$A$-tiles $U$ if and only if $f_U\in \I$.
\end{prop}

\proof We have $f_U\in \I$ if and only if
\begin{eqnarray*}
f_U(X,Y)&=&\sum_{i=1}^k\;a_iX^{\sigma_i}Y^{\tau_i}f_{T_{\lambda_i}}(X,Y)
\end{eqnarray*}
for some $a_i\in A$, $\sigma_i,\tau_i\in\R$, and
$\lambda_i\in\Lambda$.  Since
$X^{\sigma_i}Y^{\tau_i}f_{T_{\lambda_i}}(X,Y)$
is the polynomial associated to the
translation of $T_{\lambda_i}$ by the vector
$(\sigma_i,\tau_i)$, we have $f_U\in \I$ if and only
if $U=a_1\T_1+\dots+a_k\T_k$, with $\T_i$ a
translation of $T_{\lambda_i}$.  Therefore
$f_U\in \I$ if and only if
$\{T_{\lambda}:\lambda\in\Lambda\}$ $A$-tiles $U$. \qed

\begin{cor} \label{Zideal}
Let $\{T_{\lambda}:\lambda\in\Lambda\}$
be a collection of lattice tiles, let $I$ be the
ideal in $A[X^{\Z},Y^{\Z}]$ generated by the set
$\{f_{T_{\lambda}}:\lambda\in\Lambda\}$,
and let $U$ be a lattice tile such that
$f_U\in I$.  Then $\{T_{\lambda}:\lambda\in\Lambda\}$
$A$-tiles $U$.
\end{cor}

     The last result in this section shows what
happens to $f_T$ when we replace $T$ by a rescaling.

\begin{lemma} \label{scale}
Let $T$ be an $A$-weighted tile
and let $\rho$ be a positive real
number.  Then $f_{T(\rho)} =
f_T(X^{\rho},Y^{\rho})$.
\end{lemma}

\proof Let $\alpha,\beta \in \R^{\times}$.
Then $R_{\alpha\beta}(\rho)=R_{\rho\alpha,\rho\beta}$
and hence
\[f_{R_{\alpha\beta}(\rho)} \;\;=\;\;
(X^{\rho\alpha}-1)(Y^{\rho\beta}-1)
\;\;=\;\; f_{R_{\alpha\beta}}(X^{\rho},Y^{\rho}).\]
Therefore the lemma holds for $T=R_{\alpha\beta}$.
It follows from Lemma~\ref{welldef} that the lemma
holds for all tiles $T$.~\qed

\section{Tiling with rational weights}
\label{rational}
\setcounter{equation}{0}

     This section is devoted to proving
the following theorem:

\begin{theorem} \label{Q}
Let $T$ be a $\Q$-weighted tile made up of
rectangles whose corners all have rational
coordinates.  Then $T$ $\Q$-shapetiles a
square if and only if
the weighted area of $T$ is not zero.
\end{theorem}

\proof It is clear that if the weighted
area of $T$ is zero then $T$ cannot
shapetile a square with nonzero area.
Assume conversely that $T$ has nonzero
weighted area.  By rescaling and
translation we
may assume that $T$ is a lattice tile in
the first quadrant.  Let $T(\N)$ denote
the set $\{T(k):k\in\N\}$ of positive
integer rescalings of $T$.  To complete
the proof of Theorem~\ref{Q} it suffices
to prove that $T(\N)$ $\Q$-tiles a
square.  First we will prove that
$T(\N)$ $\C$-tiles a square; from this
it will follow easily that $T(\N)$
$\Q$-tiles a square.

     Since $T$ is a lattice tile in the
first quadrant, $f_T\in\Q[X,Y]$ is a polynomial
in the ordinary sense.  We begin by
interpreting the hypothesis that the
weighted area of $T$ is nonzero in terms of $f_T$.

\begin{lemma} \label{area}
There is a polynomial $\ftilde_T \in
\Q[X,Y]$ such that 
\begin{eqnarray*}
f_T(X,Y) &=& (X-1)(Y-1)\ftilde_T(X,Y).
\end{eqnarray*}
Moreover, the weighted area of $T$ is
equal to $\ftilde_T(1,1)$, and hence
$\ftilde_T(1,1) \not= 0$.
\end{lemma}

\proof Since the polynomial associated
to the unit square $S_{ij}$ is
\begin{eqnarray*}
f_{S_{ij}}(X,Y) &=& X^iY^j(X-1)(Y-1),
\end{eqnarray*}
the lemma holds for $S_{ij}$.
It follows by linearity that the lemma
holds for all lattice tiles in the first
quadrant.~\qed \m

     Let $I$ denote the ideal in
$\C[X^{\Z},Y^{\Z}]$ generated by
$\{f_{T(k)}:k\in\N\}$ and let
\begin{eqnarray*}
g_l(X,Y)&=&(X^l-1)(Y^l-1)
\end{eqnarray*}
be the polynomial associated to an $l
\times l$ square with lower left corner
$(0,0)$.
To show that $T(\N)$ $\C$-tiles a square
it suffices by Corollary~\ref{Zideal}
to show that $g_l\in I$
for some positive integer $l$.  In
order to get information about $I$ we
consider the set
$V(I)\subset\C^{\times}\times\C^{\times}$ of
common zeros of the elements of $I$.
The set $V(I)$ is essentially the union of
the lines $X=1$ and $Y=1$ with the
``shape variety'' of $T(\N)$ as defined by
Barnes \cite[\S3]{fwb1}.

     We wish to determine which points
$(\alpha,\beta)\in\C^{\times}\times\C^{\times}$ might be in
$V(I)$.  Let $m$ be the $X$-degree of $f_T$,
let $n$ be the $Y$-degree of $f_T$, and
define $\Upsilon \subset \C^{\times}$
by
\begin{eqnarray*}
\Upsilon &=& \{\zeta\in\C^{\times} : \zeta^k=1
\mbox{ for some } 1 \leq k \leq 2mn\}.
\end{eqnarray*}

\begin{lemma} \label{upsilon}
$V(I)\subset(\C^{\times}\times\Upsilon)\cup
(\Upsilon\times\C^{\times})$.
\end{lemma}

\proof Let $(\alpha,\beta)\in V(I)$, and
suppose neither $\alpha$ nor $\beta$ is
in $\Upsilon$.  By Lemma~\ref{scale}
and Lemma~\ref{area} we have
\[0 \;\;=\;\; f_{T(k)}(\alpha,\beta) \;\;=\;\;
f_T(\alpha^k,\beta^k) \;\;=\;\;
(\alpha^k-1)(\beta^k-1)
\ftilde_T(\alpha^k,\beta^k)\]
for all $k\geq1$.  Since $\alpha$
and $\beta$ aren't in $\Upsilon$ this
implies $\ftilde_T(\alpha^k,\beta^k)=0$
for $1\leq k\leq 2mn$.  Therefore by
Lemma~\ref{powers} below there exist
$c,d\in\Z$ such that $\ftilde_T(X^c,X^d)=0$.
It follows that $\ftilde_T(1,1) = 0$,
contrary to Lemma~\ref{area}.  We
conclude that if $(\alpha,\beta)\in V(I)$
then at least one of $\alpha,\beta$ must
be in $\Upsilon$.~\qed

\begin{lemma} \label{powers}
Let $K$ be a field
and let $\ftilde\in K[X,Y]$ be a nonzero
polynomial with \mbox{$X$-degree} $m-1$ and
$Y$-degree $n-1$.  Assume
there are $\alpha,\beta\in K^{\times}$
such that
\begin{enumerate}
\item $\alpha$ and $\beta$ are not $k$th
roots of 1 for any $1\leq k\leq 2mn$, and
\item $\ftilde(\alpha^k,\beta^k)=0$ for all
$1\leq k\leq 2mn$.
\end{enumerate}
Then there exist relatively prime
integers $c,d$ with $1\leq c \leq n-1$ and
$1\leq|d|\leq m-1$ such that $\ftilde(X^c,X^d)=0$.
\end{lemma}

\proof Define an $mn \times mn$
matrix $M$ whose columns are indexed by
pairs $(i,j)$ with $0\leq i\leq m-1$ and
$0\leq j\leq n-1$ by letting the $k$th entry in the
$(i,j)$ column of $M$ be $\alpha^{ik}\beta^{jk}$.
Since $\ftilde(\alpha^k,\beta^k)=0$
for $1\leq k\leq mn$, the coefficients of
$\ftilde$ give a nontrivial element of the
nullspace of $M$.  Since $M$ is essentially
a Vandermonde matrix this implies
\begin{displaymath}
0\;\;=\;\;\det(M)\;\;=\;\;
\alpha^{nm(m-1)/2}\beta^{mn(n-1)/2}\cdot\!\!\!\!\!\!
\prod_{(i,j)<(i',j')}(\alpha^{i'}\beta^{j'}-\alpha^i\beta^j)
\end{displaymath}
for an appropriate ordering of the
pairs $(i,j)$.  It follows that
$\alpha^{i'}\beta^{j'}=\alpha^i\beta^j$
for some $(i',j')\not=(i,j)$, so
$\alpha^{d_0}=\beta^{c_0}$ for some
$(c_0,d_0)\not=(0,0)$ with $|c_0|\leq n-1$ and
$|d_0|\leq m-1$.  The first assumption
implies that $c_0\not=0$ and $d_0\not=0$,
so we may assume without loss of
generality that $c_0\geq1$.

     Let $e=\gcd(c_0,d_0)$ and set
$c=c_0/e$ and $d=d_0/e$.  Then since
$(\alpha^e)^{d}=(\beta^e)^{c}$ with
$\gcd(c,d)=1$ there is a unique
$\gamma\in K$
such that $\gamma^{c}=\alpha^e$ and
$\gamma^{d}=\beta^e$.
Let $q$ be an integer such that
$1\leq q \leq 2mn/e$.  Then
by the second assumption we have
\begin{displaymath}
0\;\;=\;\;\ftilde(\alpha^{eq},\beta^{eq})\;\;=\;\;
\ftilde(\gamma^{cq},\gamma^{dq}),
\end{displaymath}
and so
$\ftilde(X^{c},X^{d}) \in K[X,X^{-1}]$ has zeros at
$X=\gamma^q$ for $1\leq q \leq 2mn/e$.
If these zeros are not distinct then for
some $1\leq r \leq 2mn/e$ we
have $\gamma^{r}=1$
and hence $1 = \gamma^{cr} = \alpha^{er}$,
which violates the first assumption.
Therefore $\ftilde(X^{c},X^{d})$ has at least
$\lfloor 2mn/e \rfloor$ distinct zeros.
On the other hand the degree of the
rational function $\ftilde(X^{c},X^{d})$ is
at most $(m-1)|c|+(n-1)|d|$, and since
$|c|=|c_0/e|\leq (n-1)/e$ and
$|d|=|d_0/e|\leq (m-1)/e$ we have
\[(m-1)|c|+(n-1)|d|\;\;\leq\;\;2(m-1)(n-1)/e\;\;<\;\;\lfloor
2mn/e \rfloor.\]
Therefore $\ftilde(X^{c},X^{d})=0$.~\qed \m

     Let $l \geq 1$ and recall that
$g_l(X,Y) = (X^l-1)(Y^l-1)$
is the polynomial associated to an $l\times l$
square with lower left corner $(0,0)$.
The set $V(g_l) \subset
\C^{\times}\times\C^{\times}$ of zeros of $g_l$
is the union of the lines $X=\zeta$ and
$Y=\zeta$ as $\zeta$ ranges over
the $l$th roots of 1.
It follows from Lemma~\ref{upsilon} that
if we choose $l$ appropriately
(say $l=(2mn)!$) then 
$V(g_l) \supset V(I)$.  This need
not imply that $g_l$ is in
$I$, but by Hilbert's
Nullstellensatz \cite[VII, Th.\,14]{zs}
we do have $g_l^k \in I$ for some
$k \geq 1$.

     To show there exists $l$
such that $g_l \in I$ we use the
theory of {\em primary decompositions}
(see, e.\,g., chapters 4 and 7 of \cite{AM}).
Let $A$ be a commutative ring with 1.
We say that the ideal $Q \subset A$ is a
{\em primary ideal} if whenever $xy \in
Q$ with $x \not\in Q$ there exists
$a \geq 1$
such that $y^a \in Q$.  By the Hilbert
basis theorem, $\C[X^{\Z},Y^{\Z}]$
is a Noetherian
ring \cite[Cor.\,7.7]{AM}.  Therefore there
are primary ideals $Q_1,\dots,Q_r$ in
$\C[X^{\Z},Y^{\Z}]$ such that $I=Q_1\cap\dots\cap Q_r$
\cite[Th.\,7.13]{AM}.  The radical ideal
\begin{displaymath}
P_i\;\;=\;\;\sqrt{Q_i}\;\;=\;\;
\{f\in\C[X^{\Z},Y^{\Z}]:f^r\in Q_i\mbox{ for some }
r\geq1\}
\end{displaymath}
of the primary ideal $Q_i$ is automatically
prime, and is called the prime
{\em associated} to $Q_i$.  We may also
characterize $P_i$ as the smallest prime
ideal containing $Q_i$.

     Since $I=Q_1\cap\dots\cap Q_r$ we
need to show that there exists $l \geq 1$
such that $g_l\in Q_i$ for all
$1\leq i\leq r$.  Observe that if $l\,|\,l'$
then $g_l\,|\,g_{l'}$.  Therefore it is
enough to show that for each $i$ there is
$l_i$ such that $g_{l_i} \in Q_i$, since in that
case we have $g_l\in I$ with
$l=\mbox{lcm}\{l_1,\dots,l_r\}$.
To accomplish this we first restrict the
possibilities for the prime ideals $P_i$.

     Let $q = (2mn)!$.  We observed above
that $g_q^k \in I$
for some positive integer $k$.  Since
$P_i \supset Q_i \supset I$ this
implies that $g_q^k \in P_i$.
Therefore some irreducible factor of
\begin{eqnarray*}
g_q(X,Y)^k &=&
\prod_{\zeta^q=1} \; (X-\zeta)^k (Y-\zeta)^k
\end{eqnarray*}
lies in the prime ideal $P_i$.  It follows that
$X-\zeta\in P_i$ or $Y-\zeta\in P_i$ for
some $\zeta \in \C^{\times}$ such that
$\zeta^q=1$.

     Assume without loss of generality
that $X-\zeta \in P_i$.  Then $P_i$
contains the prime ideal $(X-\zeta)$ generated by
the irreducible polynomial $X-\zeta$.
If $P_i \not= (X-\zeta)$ let $h$
be an element of $P_i$ which is not
in $(X-\zeta)$.  By dividing
$X-\zeta$ into $h(X,Y)$ we see that
$h(\zeta,Y) \in P_i$.  Since $P_i$ is
prime and $\C$ is algebraically closed
this implies that some linear factor $Y-\alpha$
of $h(\zeta,Y)$ is in $P_i$.
Therefore $P_i$ contains the maximal
ideal $(X-\zeta,Y-\alpha)$, so in fact
$P_i=(X-\zeta,Y-\alpha)$.  Moreover, we
must have
$\alpha \not= 0$ since $Y$ is a unit in
$\C[X^{\Z},Y^{\Z}]$.  It follows
that if $X-\zeta \in P_i$ then either
$P_i=(X-\zeta)$ or $P_i=(X-\zeta,Y-\alpha)$
for some $\alpha \in \C^{\times}$.

     We will make repeated use of the
following elementary fact about primary
ideals.

\begin{lemma} \label{primary}
Let $Q$ be a primary ideal and set $P =
\sqrt{Q}$.  If $gh \in Q$ with $h
\not\in P$ then $g \in Q$.
\end{lemma}

\proof Since $h \not\in P$ we have $h^a
\not\in Q$ for all $a \geq 1$.
Therefore by the definition of primary
ideal we have $g \in Q$.~\qed \m
     
     Assume now that $P_i=(X-\zeta)$
with $\zeta^q=1$.  Then $X^q-1$ has a
simple zero at $X=\zeta$.  Therefore
by Lemma~\ref{scale} and
Lemma~\ref{area} we have
\begin{eqnarray*}
f_{T(q)}(X,Y) &=& f_T(X^q,Y^q) \\
&=& (X^q-1)(Y^q-1)\ftilde_T(X^q,Y^q) \\
&=&(X-\zeta)h(X,Y)
\end{eqnarray*}
for some $h\in\C[X,Y]$.  Moreover we have
$h(\zeta,Y)\not=0$, since otherwise
$0=\ftilde_T(\zeta^q,Y^q) =
\ftilde_T(1,Y^q)$, which would imply
$\ftilde_T(1,1)=0$, contrary to Lemma~\ref{area}.
Therefore $h \not\in P_i = (X-\zeta)$.  It
follows by Lemma~\ref{primary} that
$X-\zeta \in Q_i$, and hence that
$g_q\in Q_i$.

     Now assume $P_i=(X-\zeta,Y-\alpha)$.
If $\alpha$
is an $r$th root of 1 for some $r \geq
1$ then $X^{qr}-1$ has a
simple zero at $X=\zeta$ and
$Y^{qr}-1$ has a simple zero at
$Y=\alpha$.  As in the previous case
this implies
\begin{eqnarray*}
f_{T(qr)}(X,Y) &=& (X^{qr}-1)(Y^{qr}-1)
\ftilde_T(X^{qr},Y^{qr}) \\
&=&(X-\zeta)(Y-\alpha)h(X,Y)
\end{eqnarray*}
for some $h\in\C[X,Y]$.  Since
$\ftilde_T(\zeta^{qr},\alpha^{qr})=\ftilde_T(1,1)\not=0$,
we have $h(\zeta,\alpha) \not= 0$, and
hence $h \not\in P_i$.
Applying Lemma~\ref{primary} we get
$(X-\zeta)(Y-\alpha)\in Q_i$, and hence
$g_{qr}\in Q_i$.  If $\alpha$ is not a
root of 1 we may choose $r \geq 1$ so that
$\ftilde_T(\zeta^{qr},\alpha^{qr})=\ftilde_T(1,\alpha^{qr})
\not=0$, since $\ftilde_T(1,1) \not= 0$
implies that $\ftilde_T(1,Y)$
has only finitely many zeros.
Then $X^{qr}-1$ has a
simple zero at $X = \zeta$ and
$Y^{qr}-1$ is nonzero at $Y = \alpha$.
By an argument similar to those used above we have
$f_{T(qr)}(X,Y)=(X-\zeta)h(X,Y)$ for
some $h\in\C[X,Y]$ such that
$h(\zeta,\alpha)\not=0$.  This implies
$h \not\in P_i$, so by
Lemma~\ref{primary} we get $X-\zeta\in Q_i$,
and hence $g_q\in Q_i$.

     We've shown now that for each
$1\leq i\leq r$ there is $l_i\geq1$ such that
$g_{l_i}\in Q_i$.  Therefore we have $g_l\in I$
with $l=\mbox{lcm}\{l_1,\dots,l_r\}$.
It follows from Corollary~\ref{Zideal} that
$T(\N)$ $\C$-tiles an $l\times l$ square.
To prove that $T(\N)$ $\Q$-tiles a square it
is sufficient to prove that $g_l$ is in the
ideal $I_0$ in $\Q[X^{\Z},Y^{\Z}]$
generated by $T(\N)$.  Equivalently, we need to
show that $g_l$ is in the $\Q$-span of the set
\begin{eqnarray*}
\E &=& \{X^i Y^j f_{T(k)} : i,j,k \in \Z, \,k
\geq 1\}.
\end{eqnarray*}
We have shown that
$g_l$ is in the $\C$-span of $\E$.
Since $g_l$ and the elements of $\E$
are all in $\Q[X^{\Z},Y^{\Z}]$, and
\begin{eqnarray*}
\C[X^{\Z},Y^{\Z}] &\cong&
\Q[X^{\Z},Y^{\Z}]\otimes_{\Q}\C, 
\end{eqnarray*}
it follows
immediately that $g_l$ is in the
$\Q$-span of $\E$.  This completes the
proof of Theorem~\ref{Q}.~\qed

\begin{cor} \label{w}
Let $T$ be a $\Z$-weighted tile made up of
rectangles whose corners all have rational
coordinates.  Assume that the
weighted area of $T$ is not zero.  Then
there exists a positive integer $w$ such
that $T(\N)$ $\Z$-tiles a square with
weight $w$.
\end{cor}

\proof By Theorem~\ref{Q} we know that
$T(\N)$ $\Q$-tiles a square $R$, so
there are rational numbers
$a_1,\dots,a_n$ and tiles
$T_1,\dots,T_n$, each a translation of
some $T(k_i) \in T(\N)$, such that $R =
a_1T_1+\dots+a_nT_n$.  Let $w \geq 1$ be
a common denominator for
$a_1,\dots,a_n$.  Then $wR =
wa_1T_1+\dots+wa_nT_n$, and $wa_i \in
\Z$ for $1 \leq i \leq n$.  Therefore
$T(\N)$ $\Z$-tiles $wR$.~\qed

\section{Tiling with integer weights}
\label{integer}
\setcounter{equation}{0}

     Let $T$ be a $\Z$-weighted lattice
tile, and assume that the weighted area
of $T$ is not zero.  By Corollary~\ref{w} we
know that $T$ $\Z$-shapetiles a
square with weight $w$ for some positive
integer $w$.  We wish to find necessary
and sufficient conditions for $T$ to
$\Z$-shapetile a square with weight
1.  To express these conditions we need
a definition.
Given $\mu\in\Q\,\cup\{\infty\}$ we say that
two lattice squares $S_{ij}$ and
$S_{i'j'}$ belong to the same {\em
$\mu$-slope class} if the line joining their
centers has slope $\mu$.  The tile $T$ can
be decomposed into a sum $T = C_1+\cdots+C_k$ of
lattice tiles
such that for each $i$ the unit lattice squares
which make up $C_i$ all belong to the same
$\mu$-slope class.

\begin{prop} \label{factor}
Let $T$ be a $\Z$-weighted lattice tile and
let $n$ be a positive integer.  Let
$c$ and $d$ be relatively prime integers and set
$\mu = -c/d$.  Then the $\mu$-slope
classes of $T$ all have weighted area divisible
by $n$ if and only if $f_T$ is an element of the
ideal $((X^d-Y^c)(X-1)(Y-1),n(X-1)(Y-1))$ in
$\Z[X^{\Z},Y^{\Z}]$.
\end{prop}

\proof The $\mu$-slope classes of $T$ all have
weighted area divisible by $n$ if and only if we
can write $T=T_1+nT_2$, where $T_1$ and $T_2$ are
$\Z$-weighted lattice tiles such
that the $\mu$-slope classes of
$T_1$ all have weighted area zero.
Write the decomposition of $T_1$ into
its $\mu$-slope classes as $T_1 = C_1+\cdots+C_k$.
Since $\mu = -c/d$ with $c$ and $d$
relatively prime, the
lattice squares $S_{ij}$ and $S_{i'j'}$
are in the same $\mu$-slope class if
and only if $S_{i'j'}$ is the translation
of $S_{ij}$ by $(dr,-cr)$ for some $r \in \Z$.
Therefore if $C_t$ is the $\mu$-slope
class of $T_1$ containing $S_{ij}$ we have
\begin{eqnarray*}
f_{C_t}(X,Y) &=& g(X^dY^{-c})X^iY^j(X-1)(Y-1)
\end{eqnarray*}
for some $g \in \Z[X^{\Z}]$.
Since the weighted area of $C_t$ is
zero we see that $0 =
\ftilde_{C_t}(1,1) = g(1)$, which implies
$X-1\,|\,g(X)$.  It follows that
$(X^dY^{-c}-1)(X-1)(Y-1)$ divides $f_{C_t}$
for $1\leq t\leq k$, and hence also that
$(X^dY^{-c}-1)(X-1)(Y-1)$ divides $f_{T_1}$.
Conversely, if $(X^dY^{-c}-1)(X-1)(Y-1)$ divides
$f_{T_1}$, it is easy to check that the $\mu$-slope
classes of $T_1$ all have weighted area zero.
It follows that the $\mu$-slope classes of $T$
all have area divisible by $n$ if and only if we
can write
\begin{eqnarray*}
f_T(X,Y) &=&
(X^dY^{-c}-1)(X-1)(Y-1)h_1(X,Y)+n(X-1)(Y-1)h_2(X,Y)
\end{eqnarray*}
for some $h_1,h_2 \in
\Z[X^{\Z},Y^{\Z}]$.  Since $Y^c$ is a unit in
$\Z[X^{\Z},Y^{\Z}]$ this is equivalent to
$f_T \in ((X^d-Y^c)(X-1)(Y-1),n(X-1)(Y-1))$.~\qed

\begin{theorem} \label{Z}
Let $T$ be a $\Z$-weighted lattice tile.
Then $T$ $\Z$-shapetiles a square if and
only if the following two conditions hold:
\begin{enumerate}
\item The weighted area of $T$ is
not zero.
\item For every $\mu\in\Q^{\times}$ the
gcd of the weighted areas of the $\mu$-slope
classes of $T$ is~1.
\end{enumerate}
\end{theorem}

\proof Let $T$ be a tile which satisfies
conditions 1 and 2.  To show that $T$
$\Z$-shapetiles a square it is sufficient
by Corollary~\ref{w} to show
that $T(\N)\cup\{wR\}$ $\Z$-tiles a square,
where $R$ is an $l\times l$ square and $l,w$
are positive integers.  Let $S=S_{00}$ be
the unit lattice
square with lower left corner $(0,0)$.  If
$T(\N)\cup\{wS\}$ $\Z$-tiles an $a
\times a$ square then by rescaling we see that
$T(\N)\cup\{wR\}$ $\Z$-tiles an $la
\times la$ square.  Therefore
it is sufficient to show that $T(\N) \cup
\{wS\}$ $\Z$-tiles a square.
Let $J$ be the ideal in
$\Z[X^{\Z},Y^{\Z}]$ generated by
$\{f_{T(k)}:k\in\N\}\cup\{w(X-1)(Y-1)\}$.
By Corollary~\ref{Zideal} it is sufficient
to show that $g_l\in J$ for some
$l \geq 1$.

     By the Hilbert basis theorem
$\Z[X^{\Z},Y^{\Z}]$ is a
Noetherian ring.  Therefore
the ideal $J$ has a primary decomposition
$J=Q_1\cap\dots\cap Q_t$.  We need to
show that there exists $l \geq 1$
such that $g_l\in Q_i$ for all $i$. 
As in the proof of Theorem~\ref{Q} it
is enough to show that for each $i$ there
is $l_i \geq 1$ such that
$g_{l_i} \in Q_i$.  Let $P_i=\sqrt{Q_i}$
be the prime associated to $Q_i$, and suppose
$w \not\in P_i$.  Then since $w(X-1)(Y-1)
\in Q_i$, by
Lemma~\ref{primary} we see that $(X-1)(Y-1)
= g_1$ is in $Q_i$.  If $w \in P_i$ then
since $P_i$ is a prime ideal it follows that
$P_i$ contains a prime integer $p$ which
divides $w$, and hence that $P_i\cap\Z=p\Z$.

     For $f\in\Z[X^{\Z},Y^{\Z}]$ let
$\fbar\in\Fp[X^{\Z},Y^{\Z}]$
be the reduction of $f$ modulo $p$,
where $\Fp=\Z/p\Z$ is the field
with $p$ elements.
Let $\Pbar_i$ be the ideal in
$\Fp[X^{\Z},Y^{\Z}]$ consisting of the reductions
modulo $p$ of the elements of $P_i$.
Since $p\in P_i$ the ideal $\Pbar_i$ is
prime.  Let $\Jbar\subset\Fp[X^{\Z},Y^{\Z}]$ be
the ideal consisting of the reductions
modulo $p$ of the elements of $J$.  Then
$\Jbar$ is generated by
$\{\fbar_{T(k)}:k\geq1\}$.  Since
$P_i\supset J$, we
have $\Pbar_i\supset\Jbar$.

     Let $K$ be an algebraic closure
of $\Fp$ and let $V(\Jbar)\subset K^{\times}\times K^{\times}$
be the set of common zeros of the
elements of $\Jbar$.
Let $m$ be the $X$-degree of $\fbar_T$,
let $n$ be the $Y$-degree of $\fbar_T$, and
define $\Upsilonbar \subset K^{\times}
\times K^{\times}$ by
\begin{eqnarray*}
\Upsilonbar &=& \{\zeta\in K^{\times}:
\zeta^k=1\mbox{ for some } 1 \leq k \leq 2mn\}.
\end{eqnarray*}

\begin{lemma} \label{subset}
$V(\Jbar) \;\subset\; (K^{\times}\times\Upsilonbar)
\cup (\Upsilonbar\times K^{\times})$.
\end{lemma}

\proof Let $(\alpha,\beta)\in V(\Jbar)$ and
suppose neither $\alpha$ nor $\beta$ is
in $\Upsilonbar$.  Then for
$1 \leq k \leq 2mn$ we have
\[0 \;\;=\;\; \fbar_{T(k)}(\alpha,\beta) \;\;=\;\;
\fbar_T(\alpha^k,\beta^k) \;\;=\;\;
(\alpha^k-1)(\beta^k-1)\ftildebarT(\alpha^k,\beta^k).\]
Since $\alpha$ and $\beta$ aren't in
$\Upsilonbar$ this implies that
$\ftildebarT(\alpha^k,\beta^k) = 0$ for
$1 \leq k \leq 2mn$.  Therefore
by Lemma~\ref{powers} there are relatively
prime integers $c,d$ with $c\geq1$ and
$d\not = 0$ such that $\ftildebarT(X^c,X^d)=0$.
Let $\A$ be the quotient ring
$\Fp[X^{\Z},Y^{\Z}]/(X^d-Y^c)$, and let
$x,y$ denote the images of $X,Y$ in $\A$.
Then $x$ and $y$ are units in $\A$
satisfying $x^d=y^c$ with $\gcd(c,d)=1$, so there
is $z=x^ay^b$ in $\A^{\times}$ such that $x=z^c$
and $y=z^d$.  Therefore the image of $\ftildebarT$
in $\A$ is given by $\ftildebarT(x,y)=
\ftildebarT(z^c,z^d)$, which equals zero since
$\ftildebarT(X^c,X^d)=0$.  It follows that
$X^d-Y^c$ divides $\ftildebarT$, and hence
that $f_T^*$ is in the ideal $(X^d-Y^c,p)$ in
$\Z[X^{\Z},Y^{\Z}]$.  Therefore
$f_T=(X-1)(Y-1)f_T^*$ is in the ideal
\[((X^d-Y^c)(X-1)(Y-1),p(X-1)(Y-1))\]
in $\Z[X^{\Z},Y^{\Z}]$.
Proposition~\ref{factor} now implies that
every $\mu$-slope class of $T$ has
area divisible by $p$.  This violates
condition 2 of the theorem, so we have a
contradiction.~\qed \m

     Set $q = (2mn)!$ and let
$V(\gbar_q) \subset K^{\times} \times K^{\times}$ be
the set of zeros of $\gbar_q$.
Since $X^q-1$
has zeros at all elements of $\Upsilonbar$,
we have $V(\gbar_q) \supset
(K^{\times}\times\Upsilonbar) \cup
(\Upsilonbar\times K^{\times})$.  Therefore
Lemma~\ref{subset} implies
$V(\gbar_q) \supset V(\Jbar)$.
Since $\Pbar_i\supset\Jbar$ we have
$V(\Jbar) \supset V(\Pbar_i)$, and
hence $V(\gbar_q) \supset
V(\Pbar_i)$.
As in Section~\ref{rational} Hilbert's
Nullstellensatz implies that $\gbar_q^k \in
\Pbar_i$ for some $k \geq 1$.  Since $\Pbar_i$
is prime and 
\begin{eqnarray*}
\gbar_q(X,Y)^k &=&
(X^q-1)^k(Y^q-1)^k
\end{eqnarray*}
we have either $X^q-1 \in \Pbar_i$
or $Y^q-1 \in \Pbar_i$.
It follows that $P_i$ contains one
of the ideals $(X^q-1,p)$ or $(Y^q-1,p)$.
We may assume without loss of generality that
$P_i\supset(X^q-1,p)$.

     By \cite[Prop.\,7.14]{AM} we
have $Q_i\supset P_i^u$ for some $u\geq1$.
Therefore it is enough to prove that for
every $u \geq 1$ there is $l \geq 1$
such that $g_l\in P_i^u$.
Let $t$ be a positive integer.
Expanding $X^{qt}-1$ in powers of
$X^q-1$ gives
\begin{eqnarray*}
X^{qt}-1&=&-1+((X^q-1)+1)^t \\
&=& \sum_{j=1}^t
\left(\!\begin{array}{c}t\\j\end{array}\!\right)
(X^q-1)^j.
\end{eqnarray*}
If we choose $t$ to be divisible by a large
power of $p$ then for small values of
$j\geq1$ the binomial coefficient
$\left(\!\begin{array}{c}t\\j\end{array}\!\right)$
is divisible by a large power of $p$.  Thus
every term in this expansion is
divisible either by a large power of $p$
or a large power of $X^q-1$.  It follows
that there
exists $t \geq 1$ such that $X^{qt}-1\in
(X^q-1,p)^u$.
Since $P_i^u\supset(X^q-1,p)^u$ we get $g_{qt} \in P_i^u$,
as required. \m

     Assume conversely that $T$ $\Z$-shapetiles a
square.  Then the weighted area of $T$ is clearly
not equal to zero, so condition 1 of Theorem~\ref{Z}
is satisfied.  We need to show that for every
$\mu \in \Q^{\times}$ the gcd of the weighted areas
of the $\mu$-slope classes of $T$ is equal to 1.
If we knew that the scale factors and
the coordinates of the translation
vectors used in shapetiling the square were all in
$\Z$, or even in $\Q$, we could prove
this using polynomials in $\Z[X^{\Z},Y^{\Z}]$.
Since we have no right to make
this assumption, we need to work in  the ring
$\Z[X^{\R},Y^{\R}]$.

     We may assume that the square
which is shapetiled by $T$ is $S=S_{00}$, the unit
square with lower left corner $(0,0)$.
We have then $S = a_1T_1+\cdots+a_kT_k$,
where $a_i \in \Z$ and each $T_i$ is a
translation of some $T(\rho_i)$.
Let $p$ be prime and suppose that
for some $\mu \in \Q^{\times}$ the
areas of the $\mu$-slope classes of
$T$ are all divisible by $p$.
Let $c,d$ be integers
such that $\gcd(c,d)=1$ and $\mu = -c/d$.
Let $\fbar_{T} \in
\Fp[X^{\Z},Y^{\Z}]$ be the
reduction of $f_{T}$ modulo $p$, and for $1\leq
i\leq n$ let $\fbar_{T_i}\in\Fp[X^{\R},Y^{\R}]$ be
the reduction of $f_{T_i}$.
Then by Proposition~\ref{factor} we see that
$(X^d-Y^c)(X-1)(Y-1)$ divides $\fbar_{T}$ (in
$\Fp[X^{\Z},Y^{\Z}]$, and hence also in
$\Fp[X^{\R},Y^{\R}]$).
Therefore by Lemma~\ref{scale} and
Lemma~\ref{translate}
we see that $\fbar_{T_i}$ is divisible by
\[(X^{\rho_i d}-Y^{\rho_i c})
(X^{\rho_i}-1)(Y^{\rho_i}-1).\]

     Define a ring homomorphism $\Psi :
\Fp[X^{\R},Y^{\R}] \ra \Fp[X^{\R}]$ by
setting $\Psi(f)=f(X^c,X^d)$.  Since
$\Psi(X^{\rho_i d}-Y^{\rho_i c}) = 0$,
the divisibility relation
from the preceding paragraph implies
that $\Psi(\fbar_{T_i})
= 0$ for $1 \leq i \leq n$.  On the
other hand, since $\fbar_S = \gbar_1 = (X-1)(Y-1)$, we have
\begin{eqnarray*}
\Psi(\fbar_S) &=& X^{c+d} -
X^{c} - X^{d} + 1,
\end{eqnarray*}
which is nonzero since $c$ and $d$
are nonzero.  Since $S =
a_1T_1+\cdots+a_kT_k$ we
have $\fbar_S =
\abar_1\fbar_{T_1}+\cdots+\abar_k\fbar_{T_k}$
with $\abar_i \in \Fp$, which gives a
contradiction.  Therefore
the areas of the
$\mu$-slope classes of $T$ can't all
be divisible by $p$, so condition 2
is satisfied.  This
completes the proof of Theorem~\ref{Z}.~\qed \m

\begin{figure}
\begin{center}
\begin{picture}(400,150)
\put(35,0){\vector(0,1){150}}
\put(15,20){\vector(1,0){150}}
\put(255,0){\vector(0,1){150}}
\put(235,20){\vector(1,0){150}}

\put(65,20){\line(0,1){30}}
\put(95,50){\line(0,1){30}}
\put(65,50){\line(1,0){30}}

\thicklines

\put(35,20){\line(0,1){30}}
\put(35,20){\line(1,0){60}}
\put(95,20){\line(0,1){30}}
\put(95,50){\line(1,0){30}}
\put(35,50){\line(1,0){30}}
\put(65,50){\line(0,1){30}}
\put(65,80){\line(1,0){60}}
\put(125,50){\line(0,1){30}}

\put(255,20){\line(0,1){110}}
\put(255,20){\line(1,0){80}}
\put(335,20){\line(0,1){50}}
\put(255,130){\line(1,0){30}}
\put(285,70){\line(1,0){50}}
\put(285,70){\line(0,1){60}}

\put(293,10){$a$}
\put(245,70){$b$}
\put(308,62){$c$}
\put(275,95){$d$}

\put(95,-10){(a)}
\put(315,-10){(b)}

\end{picture}
\end{center}
\caption{Two tiles}
\label{S}
\end{figure}
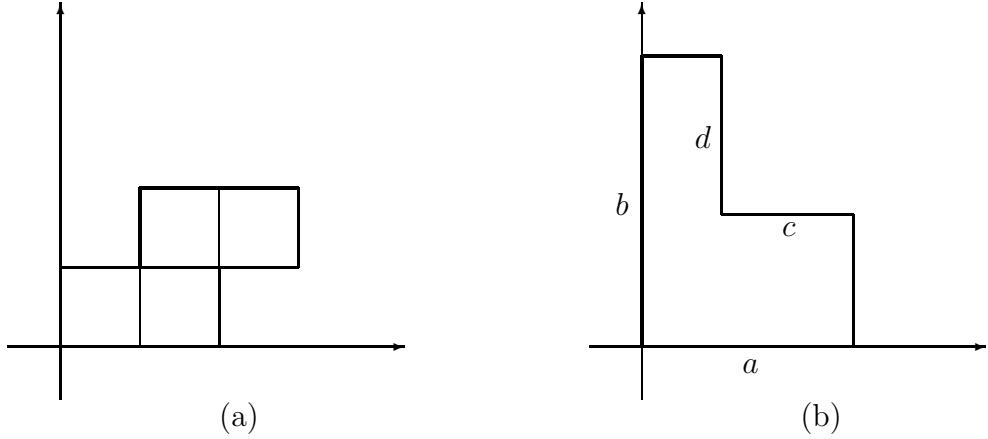

\begin{example}
Let $T$ be the lattice tile pictured in
Figure~\ref{S}a.  Since $T$ has area
$4\not=0$, it follows from Theorem~\ref{Q}
that $T$ $\Q$-shapetiles a square.
But since the nonempty 1-slope
classes of $T$ both have area 2,
Theorem~\ref{Z} implies that $T$ does not
$\Z$-shapetile a square.
\end{example}

\m
\begin{example} \label{L}
Let $a,b,c,d$ be positive integers with
$a>c$ and $b>d$.  We construct a
lattice tile $T$
by removing a $c\times d$ rectangle from
the upper right corner of an $a\times b$
rectangle, as in Figure~\ref{S}b.  The
area of $T$ is $ab-cd>0$, so the first
condition of Theorem~\ref{Z} is
satisfied.  
If $\mu>0$ there is a $\mu$-slope class
of $T$
consisting of just the upper left
corner square, while if $\mu<0$ there
is a $\mu$-slope class of $T$
consisting of just the lower left
corner square.  In either case
$T$ has a $\mu$-slope class whose
area is 1.  Therefore the
second condition of Theorem~\ref{Z} is
also satisfied, so $T$ $\Z$-shapetiles a
square.
\end{example} \m

\begin{figure}
\begin{center}
\begin{picture}(400,150)
\put(35,0){\vector(0,1){150}}
\put(15,20){\vector(1,0){150}}
\put(255,0){\vector(0,1){150}}
\put(235,20){\vector(1,0){150}}
\multiput(55,17)(20,0){6}{\line(0,1){6}}
\multiput(32,40)(0,20){6}{\line(1,0){6}}
\multiput(275,17)(20,0){6}{\line(0,1){6}}
\multiput(252,40)(0,20){6}{\line(1,0){6}}


\put(74,60){$2$}
\put(119,27){$1$}
\put(44,125){$1$}
\put(104,90){$1$}

\put(348,27){$-1$}
\put(358,47){$-1$}
\put(338,47){$-2$}
\put(338,67){$-2$}
\put(358,67){$-1$}
\put(338,87){$-1$}
\put(322,87){$0$}
\put(300,69){$-1$}
\put(318,107){$-1$}
\put(298,107){$-2$}
\put(278,107){$-2$}
\put(278,127){$-1$}
\put(298,127){$-1$}
\put(258,119){$-1$}
\put(270,39){$-1$}

\put(45,-10){Tiles with weight $+1$}
\put(265,-10){Tiles with weight $-1$}

\put(36,21){\line(1,0){120}}
\put(36,21){\line(0,1){120}}
\put(36,141){\line(1,0){60}}
\put(156,21){\line(0,1){60}}
\put(96,81){\line(1,0){60}}
\put(96,81){\line(0,1){60}}

\put(54,39){\line(1,0){80}}
\put(54,39){\line(0,1){80}}
\put(54,119){\line(1,0){40}}
\put(134,39){\line(0,1){40}}
\put(94,79){\line(1,0){40}}
\put(94,79){\line(0,1){40}}

\put(95,80){\line(1,0){40}}
\put(95,80){\line(0,1){40}}
\put(95,120){\line(1,0){20}}
\put(135,80){\line(0,1){20}}
\put(115,100){\line(1,0){20}}
\put(115,100){\line(0,1){20}}

\put(256,21){\line(1,0){80}}
\put(256,21){\line(0,1){80}}
\put(256,101){\line(1,0){40}}
\put(336,21){\line(0,1){40}}
\put(296,61){\line(1,0){40}}
\put(296,61){\line(0,1){40}}

\put(334,19){\line(1,0){40}}
\put(334,19){\line(0,1){40}}
\put(334,59){\line(1,0){20}}
\put(374,19){\line(0,1){20}}
\put(354,39){\line(1,0){20}}
\put(354,39){\line(0,1){20}}

\put(335,40){\line(1,0){40}}
\put(335,40){\line(0,1){40}}
\put(335,80){\line(1,0){20}}
\put(375,40){\line(0,1){20}}
\put(355,60){\line(1,0){20}}
\put(355,60){\line(0,1){20}}

\put(336,61){\line(1,0){40}}
\put(336,61){\line(0,1){40}}
\put(336,101){\line(1,0){20}}
\put(376,61){\line(0,1){20}}
\put(356,81){\line(1,0){20}}
\put(356,81){\line(0,1){20}}

\put(254,99){\line(1,0){40}}
\put(254,99){\line(0,1){40}}
\put(254,139){\line(1,0){20}}
\put(294,99){\line(0,1){20}}
\put(274,119){\line(1,0){20}}
\put(274,119){\line(0,1){20}}

\put(275,100){\line(1,0){40}}
\put(275,100){\line(0,1){40}}
\put(275,140){\line(1,0){20}}
\put(315,100){\line(0,1){20}}
\put(295,120){\line(1,0){20}}
\put(295,120){\line(0,1){20}}

\put(296,101){\line(1,0){40}}
\put(296,101){\line(0,1){40}}
\put(296,141){\line(1,0){20}}
\put(336,101){\line(0,1){20}}
\put(316,121){\line(1,0){20}}
\put(316,121){\line(0,1){20}}

\put(294,59){\line(1,0){40}}
\put(294,59){\line(0,1){40}}
\put(294,99){\line(1,0){20}}
\put(334,59){\line(0,1){20}}
\put(314,79){\line(1,0){20}}
\put(314,79){\line(0,1){20}}

\end{picture}
\end{center}
\caption{A $\Z$-shapetiling of a square}
\label{simple}
\end{figure}
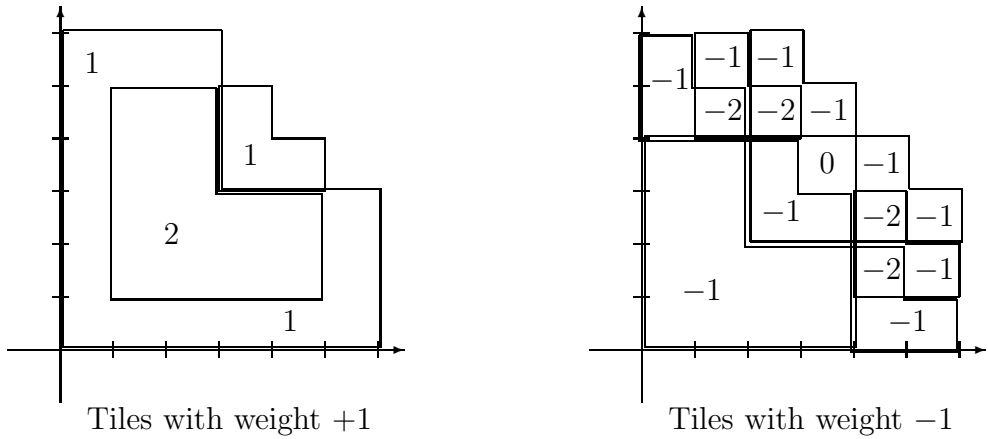

\begin{example}
The simplest case of Example~\ref{L}
occurs when $a=b=2$ and $c=d=1$.  In
this case we have $f_T(X,Y)=(1+X+Y)(X-1)(Y-1)$.  A
straightforward calculation shows that
\[\begin{array}{c}
XYg_{3}(X,Y) \,=\,
(X^3Y^3-X^2Y^2-X^4-X^4Y-X^4Y^2-Y^4-XY^4-X^2Y^4)f_T(X,Y) \\[.1cm]
+(XY-1)f_{T(2)}(X,Y)+f_{T(3)}(X,Y).
\end{array}\]
This gives the $\Z$-tiling of a $3\times3$ square
with lower left corner $(1,1)$ depicted in
Figure~\ref{simple}.  The left side of
Figure~\ref{simple} has tiles with weight 1 and
the right side has tiles with weight $-1$.  The
total weights of the tiles covering each
region are indicated.
\end{example}

\end{document}